\documentclass[12pt]{amsart}
\usepackage{etex}
\makeindex

\usepackage[ruled]{algorithm2e}
\usepackage{graphicx,tikz,url,hyperref}
\usepackage[latin1]{inputenc} \usepackage[T1]{fontenc} \usepackage{lmodern}
\usepackage{blindtext}
\usepackage{hyperref}
\usepackage{pgf}
\usepackage{enumerate}
\usepackage{epsfig,graphicx,color,amsmath,a4wide,amssymb,amsfonts,amsbsy,xy,latexsym,epsf,pstricks,verbatim,times}
\usepackage{MnSymbol}
\usepackage{pgf,color}

\usepackage{changebar}

\SetKwInOut{Input}{Input}
\SetKwInOut{Output}{Output}
\DontPrintSemicolon

\definecolor{LightGrey}{rgb}{.85,.85,.85}
\definecolor{DarkGrey}{rgb}{.5,.5,.5}

\definecolor{Blue}{rgb}{.0,.0,0.9}
\definecolor{LightBlue1}{rgb}{.2,.4,0.9}
\definecolor{LightBlue2}{rgb}{.3,.5,0.9}
\definecolor{LightBlue3}{rgb}{.4,.6,0.9}
\definecolor{LightBlue4}{rgb}{.5,.7,.9}
\definecolor{LightBlue5}{rgb}{.6,.8,.9}
\definecolor{LightBlue6}{rgb}{.7,.9,.9}

\definecolor{Red}{rgb}{.9,.0,.0}
\definecolor{LightRed1}{rgb}{0.9,.2,.4}
\definecolor{LightRed2}{rgb}{0.9,.3,.5}
\definecolor{LightRed3}{rgb}{0.9,.4,.6}
\definecolor{LightRed4}{rgb}{.9,.5,.7}
\definecolor{LightRed5}{rgb}{.9,.6,.8}
\definecolor{LightRed6}{rgb}{.9,.7,.9}

\xyoption{all}
\usepackage[english]{babel}

\newcounter{noalgo}[section]
\newdimen\indentalgo
\newdimen\indentalgodec\indentalgo=0.0mm\indentalgodec=10mm

\def\<<{\leavevmode
  \raise0.28ex\hbox{$\scriptscriptstyle\langle\!\langle$}\nobreak
  \hskip -.6pt plus.3pt minus.2pt\,}
\def\>>{\,\nobreak\hskip -.6pt plus.3pt minus.2pt
  \raise0.28ex\hbox{$\scriptscriptstyle\rangle\!\rangle$}}

\newtheorem{exple}{Example}

\newtheorem{theorem}{Theorem}

\providecommand{\myproofname}{Proof}

\usepackage{enumitem}

\usepackage{scalerel,stackengine}
\stackMath
\newcommand\widecheck[1]{%
\savestack{\tmpbox}{\stretchto{%
  \scaleto{%
    \scalerel*[\widthof{\ensuremath{#1}}]{\kern-.6pt\bigwedge\kern-.6pt}%
    {\rule[-\textheight/2]{1ex}{\textheight}}
  }{\textheight}%
}{0.5ex}}%
\stackon[1pt]{#1}{\scalebox{-1}{\tmpbox}}%
}

\begin{document}

\title{Rational Points of some genus 3 curves from the rank $0$ quotient strategy}

\author{Tony Ezome}
\address{Tony Ezome, \'Ecole Normale Sup\'erieure,
D{\'e}partement de math{\'e}matiques,
BP 17 009 Libreville, Gabon.}
\email{tony.ezome@gmail.com}

\author{Brice Miayoka Moussolo}
\address{Brice Miayoka Moussolo, Universit{\'e} Marien Ngouabi, 
Facult\'e des Science et Techniques BP 69, Brazzaville, Congo}
\email{bricemiayo@gmail.com}

\author{R\'egis Freguin Babindamana}
\address{R\'egis Freguin Babindamana, Universit{\'e} Marien Ngouabi, 
Facult\'e des Science et Techniques BP 69, Brazzaville, Congo}
\email{regis.babindamana@yahoo.fr}

\date{\today}

\maketitle
\setcounter{tocdepth}{2} 

\begin{abstract}
In 1922, Mordell conjectured
that the set of rational points on a smooth
curve $C$ over $\mathbb{Q}$ with genus $g\ge 2$ is finite.
This has been proved by Faltings in 1983.
However, Coleman determined in 1985 an upper bound of
$\#C(\mathbb{Q})$ by following Chabauty's approach
which considers the special case when the Jacobian
variety of $C$ has Mordell-Weil rank $< g$.
In 2006,  Stoll improved the Coleman's bound.
Balakrishnan with her co-authors in \cite{Balakrishnan_et_al} implemented
the Chabauty-Coleman method to
compute the rational points of genus $3$ hyperelliptic curves.
Then, Hashimoto and Morrison \cite{Hashimoto-Morrison}
did the same work for Picard curves.
But it happens that this work has not yet been done
for all genus $3$ curves.
In this paper, we describe an algorithm
to compute the complete set of rational
points $C(\mathbb{Q})$ for any genus $3$
curve $C/\mathbb{Q}$ that is a degree-$2$ cover of a
genus 1 curve whose Jacobian has rank $0$.
We implemented this algorithm in Magma, and we ran 
it on approximately $40,000$ curves selected from
databases of plane quartics and genus $3$ hyperellitic curves.
We discuss some interesting examples, and
we exhibit curves for which the number of rational
points meets the Stoll's bound.
\end{abstract}



\section{Introduction}

Given a projective, smooth, absolutely integral curve $C$ over
$\mathbb{Q}$, we are interested in determining the set
$C(\mathbb{Q})$ of rational points on $C$.
The genus $g$ of $C$, which is a nonnegative
integer depending on $C$ up to birational equivalence,
is an important data. If $g=0$, then either $C(\mathbb{Q})=\emptyset$ or
$C$ is isomorphic to the projective line.
If $g=1$ and $C(\mathbb{Q})\neq\emptyset$,
then $C$ is an elliptic curve.
In the latter case, a famous theorem by Mordell in \cite{Mordell} certifies that
$C(\mathbb{Q})$ is 
a finitely generated abelian group. This means that $C(\mathbb{Q})$
can be described only from a finite number of its points.
Moreover, at the end of his paper Mordell conjectured
that if $g$ is greater than or equal to $2$, then
$C(\mathbb{Q})$ is finite.
In 1929, Weil \cite{Weil} generalized the Mordell's theorem to 
all abelian varieties over number fields.
And then, Faltings \cite{Faltings} proved the Mordell's conjecture in 1983.
But Falting's proof is not effective.
Actually, the problem of constructing an algorithm which computes
the rational points of a given curve with genus $\ge 2$ is of
topical interest. There are some methods adapted to special
families of curves, but the problem is difficult in general.
One of the main known methods is based on the work by Chabauty \cite{Chabauty}
and Coleman \cite{Coleman_Duke}.
Assuming that $p$ is a prime of good reduction for
$C$, we denote by $P$ a point in $C(\mathbb{Q}_p)$. The embedding
$$ \xymatrix{
\iota : C  \ \   \ar@{^(-{>}}[r] &
 J\\
    Q  \ar@{|-{>}}[r] & [Q-P]
}$$
induces an isomorphism 
$\iota^* : H^0(J_{\mathbb{Q}_p},\Omega^1)
\longrightarrow H^0(C_{\mathbb{Q}_p},\Omega^1)$
between spaces of holomorphic differential forms
on $C$ and its Jacobian variety $J$.
Denote by $\overline{J(\mathbb{Q})}$ the $p$-adic closure 
of $J(\mathbb{Q})$ in $J(\mathbb{Q}_p)$.
Chabauty \cite{Chabauty} proved that if the rank of
the Jacobian of $C$ is less than $g$, then the intersection 
$\iota (C(\mathbb{Q}_p))\cap \overline{J(\mathbb{Q})}$
is a finite set. And then, Coleman \cite{Coleman_Duke}
proposed an effective version of this result by
using his theory of $p$-adic integration on curves.
Indeed, from properties of the pairing
$$ \xymatrix{
\langle, \rangle : H^0(C_{\mathbb{Q}_p},\Omega^1) \times
J(\mathbb{Q}_p)  \ar@{-{>}}[r] & 
\mathbb{Q}_p\\
    (\omega,[\sum P_i-Q_i])  \ar@{|-{>}}[r] & \sum\int_{P_i}^{Q_i} \omega,
}$$
 and
the fact that $H^0(C_{\mathbb{Q}_p},\Omega^1)$ is $g$-dimensional
while $J$ has Mordell-Weil rank $<g$,
one deduces that there exist
linearly independant differentials $\omega_1,\ldots,\omega_{k}
\in H^0(C_{\mathbb{Q}_p},\Omega^1)$ such that $\langle\omega_i,[D]\rangle=0$ for all
divisor class $[D]$ in $J(\mathbb{Q})$.  
Coleman actually computed all the $\mathbb{Q}$-rational points of $C$
by investigating integrals of the $\omega_i$'s on
$J(\mathbb{Q}_p)$ which vanish on $J(\mathbb{Q})$.
That is why these $1$-forms are usually called the
\textit{annihilating differentials}. Furthermore, he obtained
an upper bound of the number of rational points.

\begin{theorem}[\cite{Coleman_Duke}]\label{thm:1}
Let $C$ be a smooth curve over $\mathbb{Q}$ of genus $g \ge 2$
whose Jacobian has Mordell-Weil rank $r<g$, and $p$ is a prime
of good reduction such that $p>2g$. Denote
by $\overline{C}$ the reduction of $C$ modulo $p$. Then 
\begin{equation}\label{eq:1}
\#C(\mathbb{Q})\le \# \overline{C} (\mathbb{F}_{\!p})+2g-2.
\end{equation}
\end{theorem}

Later on Stoll \cite{Stoll06} improved this bound.

\begin{theorem}[\cite{Stoll06}]\label{thm:2}
Let $C$ be a smooth curve over $\mathbb{Q}$ of genus $g \ge 2$
whose Jacobian has Mordell-Weil rank $r<g-1$, and $p$ is a prime
of good reduction such that $p>2r+2$. Denote
by $\overline{C}$ the reduction of $C$ modulo $p$. Then 
\begin{equation}\label{eq:20}
\#C(\mathbb{Q})\le \# \overline{C} (\mathbb{F}_{\!p})+2r.
\end{equation}
\end{theorem}

Balakrishnan with her co-authors in \cite{Balakrishnan_et_al} implemented
the Chabauty-Coleman method to
compute the rational points of genus $3$ hyperelliptic curves
selected from the database \cite{Database_Sutherland1}.
On the other hand, Hashimoto and Morrison \cite{Hashimoto-Morrison}
did the same work for Picard curves
selected from the database \cite{Database_Sutherland2}.
But it happens that this work has not yet been done
for all genus $3$ curves. In this paper, we describe an algorithm
to compute the complete set of rational
points $C(\mathbb{Q})$ for any genus $3$
curve $C/\mathbb{Q}$ that is a degree-$2$ cover of a
genus 1 curve $D/\mathbb{Q}$ whose Jacobian has rank $0$.
The first step is
to check whether $D$ possesses a $\mathbb{Q}$-rational point or not.
When $D(\mathbb{Q})=\emptyset$, obviously $C(\mathbb{Q})=\emptyset$.
If we find a rational point $P_0$,
then $(D,P_0)$ is an elliptic curve.
Mordell's theorem tells us that $D(\mathbb{Q})$
is a finite group. Furthermore, Mazur [\cite{Mazur77}, \cite{Mazur78}] proved
that this group is isomorphic to one of the following fifteen groups:
$$
\mathbb{Z}/n\mathbb{Z} \quad \text{ with } \quad 1\le n\le 10 \ \text{ or } \ n=12,
$$
and 
$$
\mathbb{Z}/2\mathbb{Z} \times \mathbb{Z}/2n\mathbb{Z} \quad
\text{ with } 1\le n\le 4.
$$
One may use the Magma's function $\textit{MordellWeilGroup}$
to decide the exact nature of $D(\mathbb{Q})$.
Since the quotient map $\psi : C \longrightarrow D$ 
is defined over $\mathbb{Q}$, we have
$\psi(C(\mathbb{Q}))\subseteq D(\mathbb{Q})$ and
one easily deduces $C(\mathbb{Q})$ from $\psi^{-1}(D(\mathbb{Q}))$.
Note that the rank $0$ quotient strategy is not specific
to genus $3$ curves. Indeed, it is easily seen that
the algorithm described in
Section \ref{sect:algo} can be used to
compute the rational points of any curve of genus $\ge 2$
satisfying the requirements on the inputs.
For instance, a specific genus 2 case is discussed by
Siksek in \cite{siksek}.
In addition, the rank zero quotient strategy 
could be combined with the Chabauty-Coleman method
in order to implement an efficient Point-Counting algorithm
which takes advantage of strenghts of each of both procedures so that
if the implemented function fails to find the annihilating
differentials as required in the Chabauty-Coleman method, it could switch
in searching symmetries and elliptic quotients. It seems that this trick
has been used for the implementation of the Magma's Chabauty function
concerning genus $2$ curves.

Our paper is organized as follows. After 
this introduction, we recall some important properties of 
genus $3$ curves in Section \ref{sect:loci.ternary}.
Then, we describe the rank $0$ quotient strategy in Section \ref{sect:algo}.
We discuss our implementation in Section \ref{sect:implem}.
Section \ref{section:1} is devoted to illustrations.
We present some interesting
examples, and we exhibit curves for which the number of rational
points meets the Stoll's bound.

\subsection*{Acknowledgments}
This study has been carried out with financial support from the French State,
managed by CNRS in the frame of the \textit{Dispositif de Soutien aux
Collaborations avec l'Afrique subsaharienne}
(via the REDGATE Project and the IRN AFRIMath).
Experiments presented in this paper were carried
out by using the Magma
software through the account of the second author at 
Boston University Library.  The first two authors were
supported by Simons Foundation via the PREMA project.

We are very grateful to Jennifer Balakrishnan
for many helpful discussions and for having brought our attention to
recent papers on Coleman integration on curves
and Point-Counting algorithms. We thank
Oana Padurariu for helpfull suggestions.
We are also very grateful to Steffen M\"uller for
his comments on early versions of this work.
We would like to thank Andrew Sutherland 
for his comments on this work and for 
sharing with us some of the data he used
with his co-authors in \cite{Sutherland_Alii}. 

\section{The rank $0$ quotient strategy}\label{sect:loci.ternary}

As mentioned in the introduction,
this algorithm can be used to
compute the rational points of any curve of genus $\ge 2$
satisfying the requirements on the inputs.
In this section, we first recall the relevant properties
of genus $3$ curves, and then we present the algortihm.
We describe our implementation and experiments
using curves selected from the database \cite{Database_Sutherland1}
and from another database constructed by the
authors of \cite{Sutherland_Alii}.

\subsection{Genus $3$ curves are elligible}\label{sect:loci.ternary}

While families of curves with genus $g\le 2$ are homogenous
in the sense that two curves with the same genus over
an algebraically closed field are of the same nature
(rational, elliptic, hyperelliptic), things are different when
the genus is greater than $2$.
For instance, a (projective, smooth, absolutely integral) curve $C$
of genus $3$ may be either hyperelliptic or
nonhyperelliptic, $i.e$, a plane quartic.
Assume that this curve $C$
is a degree-$2$ cover of a genus 1 curve $D$.
It is known that the Jacobian variety of $C$
is isogenous to the product of the Jacobian variety of
$D$ by an abelian surface, see
for instance \cite{Ritzenthaler_Romagny}.
When $C$ is also a hyperelliptic curve, it is defined by
$C : y^2=f(x^2)$, where $f$ is a polynomial of degree
$4$ and $D$ is the locus $y^2=f(x)$.
And if $C$ is a plane quartic, then it is the locus of
a ternary quartic of the form $Y^4-h(X,Z)Y^2+r(X,Z)$
so that $D$ is given by $D : Y^2-h(X,Z)Y+r(X,Z)=0$.
The present paper actually describes how to compute
the rational points of $C$ from those of $D$ in the case when 
the Jacobian of $D$ has Mordell-Weil rank $0$. This is quite suitable
for genus $3$ hyperelliptic curves endowed with an extra involution
giving rise to a rank $0$ quotient.
Ciani quartics whose Jacobians have rank $\le 2$ are also elligible.
Recall that a Ciani quartic is a smooth curve defined by an equation
of the form
$$a_1X^4 + a_2 Y^4 + a_3 Z^4 + 2(b_1 Y^2Z^2
+ b_2X^2Z^2 + b_3 X^2Y^2)=0.$$
It is known that the Jacobian variety of such a curve
is isogenous to the product
of three quotient elliptic curves,
see \cite[Section 2.3]{Ritzenthaler_Lachaud}.
Therefore, if the Jacobian has rank $\le 2$, then
at least one of these elliptic curves has rank $0$.

\subsection{The algorithm}\label{sect:algo}

A simple presentation of the rank $0$ quotient strategy
is as follows.

\begin{algorithm}\label{algo:5}
 \caption{  \textsc{Rank $0$ quotient strategy}}

  \Input{ A curve $C/\mathbb{Q}$  and an involution $\sigma$ on $C$
  defined over $\mathbb{Q}$ such that the quotient
  $D:=C/\langle\sigma\rangle$ has 
  genus $1$ and its Jacobian has rank $0$.}

\vspace{.099cm}
               
  \Output{The set $S$ of $\mathbb{Q}$-rational points of $C$.
  }

\vspace{.3cm}  

Set $\psi : C \longrightarrow D$ the quotient map, and
compute $D(\mathbb{Q})$.\\
If $D(\mathbb{Q})= \emptyset$, then $S=\emptyset$.
Otherwise, compute $\psi^{-1}(P)$ for every $P \in D(\mathbb{Q})$.\\
Set $\Omega(P):=\psi^{-1}(P) \cap C(\mathbb{Q})$,
then $S=\bigcup_{P\in D(\mathbb{Q})}\Omega(P)$.

\vspace{.3cm}

 \Return{$S$}

\end{algorithm}

\subsection{Implementation}\label{sect:implem}

We implemented this algorithm
in Magma. The code is available
on the GitHub repositories of the second author\\
$\mathrm{\href{https://github.com/Brice202145/Strategie_Rang_quotient}{https://github.com/Brice202145/Strategie\_Rang\_quotient}}$.\\

\begin{figure}\caption{Statistics obtained when running Algorithm \ref{algo:5}
on elligible \textbf{38,564} quartic covers selected from a database
constructed in \cite{Sutherland_Alii}.}\label{fig:1}
\begin{center}
\includegraphics[height=8.94cm]{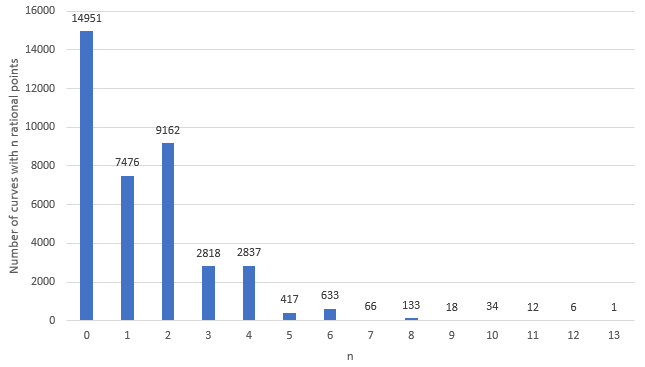} 
\end{center}
\end{figure}

We ran our code on $38,564$ plane
quartics selected from a database
constructed by Sutherland and his
co-authors in \cite{Sutherland_Alii}.
Actually, our original dataset was made of $100,000$
plane quartics which had to pass identification and selection
processes. We first tested whether their equations show certain symmetries
($i.e$ if these equations are invariant when
replacing $x$ by $-x$, $y$ by $-y$, or $z$ by $-z$)
in order to identify
the good involution for each of the curves.
In fact, we started by applying linear transformations
to these equations in order to ease the search for symmetries.
When a symmetry was found,
we checked whether the associated quotient curve has genus $1$.
We got the confirmation that all the $100,000$ curves
were symmetric with genus $1$ quotients. The selection process
consisted in computing the rank of the Jacobians of these quotient curves.
It turns out that only $38,564$ curves have a rank $0$ quotient curve,
among them there are $17,404$ Ciani quartics.
Running Algorithm \ref{algo:5} reveals that most of these curves
have no rational points.
The statistics resulting from the selection process
are summarized in table \ref{table:1}. In addition,
the statistics obtained when running Algorithm \ref{algo:5}
on the elligible $38,564$ quartic covers
are presented in Figure \ref{fig:1}.

\begin{figure}
\caption{Statistics obtained when running Algorithm \ref{algo:5}
on \textbf{130} hyperelliptic curves selected from the database
\cite{Database_Sutherland1}.}\label{fig:2}
\begin{center}
\includegraphics[height=8.94cm]{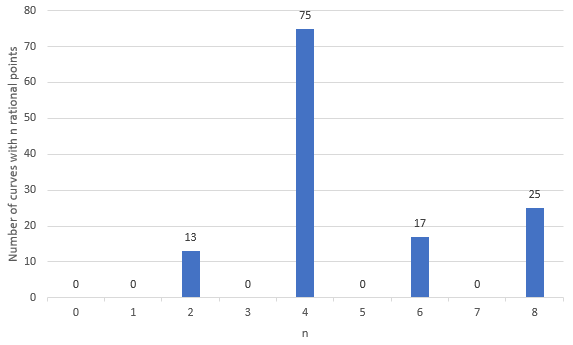} 
\end{center}
\end{figure}

We also ran our code on hyperelliptic curves
selected from the database
\cite{Database_Sutherland1}. Actually, we
applied the identication and selection processes
that we described previously to all the
$67, 879$ hyperelliptic curves of this database. It turns out that
only $130$ have a genus $1$ rank $0$ quotient curve,
and each of them possesses at least one rational point.
The statistics obtained when running Algorithm \ref{algo:5}
on the elligible $130$ hyperelliptic curves
are presented in \ref{fig:2}.
All the codes as well as the results obtains
during our experiments about
planes quartics and hyperelliptic curves are available
on the GitHub repositories of the second author\\
$\mathrm{\href{https://github.com/Brice202145/Strategie_Rang_quotient}{https://github.com/Brice202145/Strategie\textunderscore\!\!Rang\textunderscore\!\!quotient}}$.\\


\begin{table}
\caption{Statistics resulting from
our selection procedure on $100,000$ plane quartics
selected from a
database constructed in \cite{Sutherland_Alii}.}\label{table:1}
\begin{center}
$\begin{array}{|c|c|c|}\hline
{\begin{array}{l}
\text{Curves \ having\ a \ genus \ $1$}\\
\text{rank \ $0$ \ quotient}
\end{array}} & 
{\begin{array}{l}
\text{Numbers \ of \ curves}\\
\text{having \ a \ genus \ } 1\\
\text{rank $1$ \ quotient}
\end{array}} &
{\begin{array}{l}
\text{Numbers \ of \ curves}\\
\text{having \ a \ genus \ } 1\\
\text{ quotient\ of \ rank $\ge 2$}
\end{array}}\\ \hline
{\begin{array}{c|c}
{\begin{array}{c}
{\text{Numbers \ of}}\\
{\text{curves \ whose}}\\
{\text{quotient \ has \ no}}\\
{\text{rational \ points:}}\\
\\
{\mathbf{7534}} \end{array}} &
{\begin{array}{c}
{\text{Numbers \ of}}\\
{\text{curves \ whose}}\\
{\text{quotient \ is \ an}}\\
{\text{elliptic \ curve:}}\\
\\
{{\bf \Large 31030}} \end{array}}
\end{array}} &
{\mathbf{57347}} & {\mathbf{4089}}\\ \hline
\end{array}$
\end{center}
\end{table}

\vspace{.3cm}

\section{Examples}\label{section:1}

Theorem \ref{thm:1} and
Theorem \ref{thm:2} give upper bounds of
the number of rational points
on smooth curves of genus $g\ge 2$
under certain conditions.
In practice, when computing the rational points of these
curves one observes that some of them have 
a number of rational points which meets these bounds,
but some others do not.
Such observations have already been illustrated in 
\cite{Gajovic}. The following examples
also illustrate this point.

\begin{exple}[Rational points on Fermat curves]\label{exple:1}
Let $k$ be a positive integer, and let
$a_1, a_2, a_3$ nonzero rational numbers which are
$4k$-powers. By Fermat's last theorem,
there are no triple $(X,Y,Z) \in \mathbb{Q}^3$ such that 
$XYZ\ne 0$ and $(X,Y,Z)$ is a root of the trinomial
$a_1X^{4k}+a_2Y^{4k}-a_3Z^{4k}$.
Hence, the only $\mathbb{Q}$-rational points on the curve
$$C_k : X^{4k}+Y^{4k}-Z^{4k}=0$$
are the trivial ones: $(0,1,1)$, $(0,1,-1)$,
$(1,0,1),$ and $(1,0,-1)$. 
We conclude that the number of rational points on
$C_1$ meets the Stoll's bound for all prime
$p>2$. But $\#C_k(\mathbb{Q})$ never meets the Coleman's bound
for any $k$.
\end{exple}

\begin{exple}[Rational points of a Ciani quartic]\label{exple:2}
Consider the Ciani quartic
$$
C : 15X^4 + 112Y^4 +112Z^4 + 288Y^2 Z^2 - 88X^2Z^2 - 88X^2Y^2= 0
$$
Let 
$$\sigma_1: (X,Y, Z) \mapsto (-X, Y, Z), \quad
\sigma_2 : (X,Y, Z) \mapsto (X, -Y, Z), \ \text{ and } \
\sigma_3 : (X,Y, Z) \mapsto (X, Y, -Z)$$ be the three
canonical involutions of $C$.
The associated quotient elliptic curves are
$$E_1 : 112X^4 + 112Y^4 + 15Y^2Z^2 - 88ZY^3 - 88YZX^2 + 288X^2Y^2 = 0,$$
$$E_2 : 112Y^4 + 112Z^4 + 15X^2 Z^2 - 88XZ^3 - 88XZY^2 + 288Z^2 Y^2 = 0,$$
$$E_3 : 112X^4 + 112Z^4 + 15Y^2 X^2 - 88YX^3 - 88YXZ^2 + 288X^2Z^2 = 0.$$
Set $E : y^2 = x^3+\frac{7}{2}x^2- \frac{15}{16} x$. By using the following transformations
$$E_1 \longrightarrow E : (X,Y,Z)
\longmapsto \left(\frac{28}{15}X^2Y + \frac{58}{15}Y^3 - Y^2Z,
\frac{56}{15}X^3 +\frac{116}{15}XY^2 - 2XYZ, -\frac{8}{15}Y^3\right),$$
$$E_2 \longrightarrow E : (X,Y,Z)
\longmapsto \left(\frac{28}{15}Y^2Z + \frac{58}{15} Z^3 - Z^2X, \frac{56}{15}Y^3
+ \frac{116}{15} YZ^2- 2XYZ, -\frac{8}{15} Z^3\right),$$
$$
E_3 \longrightarrow E : (X,Y,Z)
\longmapsto \left(\frac{28}{15}Z^2X + \frac{58}{15} X^3 - X^2Y, \frac{56}{15}Z^3
+ \frac{116}{15} ZX^2- 2XYZ, -\frac{8}{15} X^3 \right),$$
we see that $E_1$, $E_2$ and $E_3$ are isomorphic to $E$.
Magma tells us that $E$ has rank $0$.
By running our code, we obtained
$$
E(\mathbb{Q})=
\left\{ \infty, (0, 0), \left(\frac{3}{4}, -\frac{3}{2}\right),
\left(-\frac{3}{4}, \frac{3}{2} \right),
\left(\frac{1}{4}, 0\right), \left(\frac{5}{4}, -\frac{5}{2}\right),
\left(\frac{5}{4}, \frac{5}{2}\right), \left(-\frac{15}{4}, 0\right)\right\},
$$
and
$$
C(\mathbb{Q})=
\left\{(2:1: 0), (-2:1:0), (2:0:1), (-2:0:1)\right\}.
$$
Hence, the Stoll's bound is sharp at $17$, since $C$ has rank $0$ and
its reduction modulo $17$ satisfies
$$
\overline{C}(\mathbf{F}_{\!17})=
\left\{(2 : 1 : 0),  (15 : 1 : 0), (2 : 0 : 1), (15 : 0 : 1)\right\}.
$$
\end{exple}

\begin{exple}[Rational points of a genus $3$ hyperellitic curve]\label{exple:4}
The hyperelliptic curve 
$$ C : Y^2Z^6 = X^8 + 2X^4 Z^4 - 4X^2 Z^6 + Z^8$$
is invariant under the symmetry $\sigma : (X,Y,Z) \mapsto (-X,Y,Z)$,
and the quotient $D:=C/\langle \sigma \rangle$
is a genus $1$ rank $0$ curve. By running our code, we found
$$
D(\mathbb{Q})=
\left\{\infty, (0 : \frac{1}{2} : 1), (0 : -\frac{1}{2} : 1),
(1 : \frac{1}{2} : 1), (1 : -\frac{1}{2} : 1)\right\},
$$
and
$$
C(\mathbb{Q})=
\left\{ (1 : -1 : 0), (1 : 1 : 0), (-1 : 0 : 1),
(0 : -1 : 1), (0 : 1 : 1), (1 : 0 : 1)\right\}.
$$
\end{exple}

\begin{exple}[The case of a genus $3$
nonhyperelliptic curve which is not a Ciani quartic]\label{exple:3}
Consider the plane quartic
$$ C : 3X^4 - 2X^2 Y^2 - Y^4 + 4X^3 Z - 4XY^2 Z + 2X^2 Z^2
- 2Y^2 Z^2 + 2XZ^3 + Z^4.$$
This curve is invariant under the involution $\sigma_2
: (X,Y,Z)\mapsto (X,-Y,Z)$, but not under
$\sigma_1
: (X,Y,Z)\mapsto (-X,Y,Z)$ nor $\sigma_3
: (X,Y,Z)\mapsto (X,Y,-Z)$.
The quotient $E:=C/\langle \sigma_2 \rangle$
is an elliptic curve of rank $0$. By running our code, we found
$$
E(\mathbb{Q})=
\left\{\infty, (-1: 0 :1),
\left(-\frac{1}{2}: \frac{3}{4} : 1 \right), \left(\frac{1}{2} :
\frac{1}{4}:1\right)\right\},
$$
and
$$
C(\mathbb{Q})=
\left\{(-1:0: 1), (-1: 1 :0), (1:1:0),
\left(-\frac{1}{2}: -\frac{1}{2} : 1 \right), \left(-\frac{1}{2} :
\frac{1}{2}:1\right)\right\}.
$$
We conclude that 
$\#C(\mathbb{Q})$ cannot meet the Coleman's bound,
since the reductions of $C$ satisfy
$\#\overline{C}(\mathbf{F}_{p})\ge 2$ for any prime $p$.
\end{exple}

\bibliographystyle{plain}

\bibliography{Quotient-Strategy}

\begin{thebibliography}{10}

\bibitem{Balakrishnan_et_al}
Jennifer~S. {Balakrishnan}, Francesca {Bianchi}, Victoria {Cantoral-Farf\'an},
  Mirela {\c{C}iperiani}, and Anastassia {Etropolski}.
\newblock {Chabauty-Coleman experiments for genus 3 hyperelliptic curves}.
\newblock In {\em Research directions in number theory. Women in numbers IV.
  Proceedings of the women in numbers, WIN4 workshop. Banff International
  Research Station, Banff, Alberta, Canada, August 14--18, 2017}, pages 67--90.
  Cham: Springer, 2019.

\bibitem{Database_Sutherland1}
A.~{Booker}, D.~{Platt}, and A.~{Sutherland}.
\newblock A database of nonhyperelliptic genus 3 curves over q.
\newblock Available at \begin{verbatim}
  http://math.mit.edu/~drew/gce_genus3_hyperelliptic.txt \end{verbatim}.

\bibitem{Chabauty}
Claude {Chabauty}.
\newblock {Sur les points rationnels des courbes alg\'ebriques de genre
  sup\'erieur \`a l'unit\'e}.
\newblock {\em {C. R. Acad. Sci., Paris}}, 212:882--885, 1941.

\bibitem{Coleman_Duke}
Robert~F. {Coleman}.
\newblock {Effective Chabauty}.
\newblock {\em {Duke Math. J.}}, 52:765--770, 1985.

\bibitem{Faltings}
G.~Faltings.
\newblock Endlichkeitss\"atze f\"ur abelsche variet\"aten \"uber zahlk\"orpern.
\newblock {\em Invent. Math.}, 3(73):349--366, 1983.
\newblock An optional note.

\bibitem{Sutherland_Alii}
Francesc Fit{\'e}, Kiran~S. Kedlaya, and Andrew~V. Sutherland.
\newblock Sato-{Tate} groups of abelian threefolds: a preview of the
  classification.
\newblock In {\em Arithmetic, geometry, cryptography and coding theory, AGC2T,
  17th international conference, Centre International de Rencontres
  Math\'ematiques, Marseilles, France, June 10--14, 2019}, pages 103--129.
  Providence, RI: American Mathematical Society (AMS), 2021.

\bibitem{Gajovic}
A.~{Gajovi\'c}.
\newblock Curves with sharp chabauty-coleman bound.
\newblock To appear.
\newblock The proceedings volume for the Simons Collaboration "Arithmetic
  Geometry, Number Theory, and Computation".

\bibitem{Hashimoto-Morrison}
S.~{Hashimoto} and T.~{Morrison}.
\newblock Chabauty-coleman computations on rank 1 picard curves.
\newblock Available at https://arxiv.org/pdf/2002.03291.pdf.

\bibitem{Ritzenthaler_Lachaud}
Gilles {Lachaud} and Christophe {Ritzenthaler}.
\newblock {On some questions of Serre on abelian threefolds}.
\newblock In {\em Algebraic geometry and its applications. Dedicated to Gilles
  Lachaud on his 60th birthday. Proceedings of the first SAGA conference,
  Papeete, France, May 7--11, 2007}, pages 88--115. Hackensack, NJ: World
  Scientific, 2008.

\bibitem{Mazur77}
B.~Mazur.
\newblock Modular curves and the {Eisenstein} ideal.
\newblock {\em Publ. Math., Inst. Hautes {\'E}tud. Sci.}, 47:33--186, 1977.

\bibitem{Mazur78}
B.~Mazur.
\newblock Rational isogenies of prime degree. ({With} an appendix by {D}.
  {Goldfeld}).
\newblock {\em Invent. Math.}, 44:129--162, 1978.

\bibitem{Mordell}
L.~J. {Mordell}.
\newblock {On the rational solutions of the indeterminate equations of the
  third and fourth degrees.}
\newblock {\em {Proc. Camb. Philos. Soc.}}, 21:179--192, 1922.

\bibitem{Ritzenthaler_Romagny}
Christophe Ritzenthaler and Matthieu Romagny.
\newblock On the {Prym} variety of genus 3 covers of genus 1 curves.
\newblock {\em {\'E}pijournal de G{\'e}om. Alg{\'e}br., EPIGA}, 2:8, 2018.
\newblock Id/No 2.

\bibitem{siksek}
Samir Siksek.
\newblock Chabauty and the {Mordell}-{Weil} sieve.
\newblock In {\em Advances on superelliptic curves and their applications.
  Based on the NATO Advanced Study Institute (ASI), Ohrid, Macedonia, 2014},
  pages 194--224. Amsterdam: IOS Press, 2015.

\bibitem{Stoll06}
Michael Stoll.
\newblock Independence of rational points on twists of a given curve.
\newblock {\em Compos. Math.}, 142(5):1201--1214, 2006.

\bibitem{Database_Sutherland2}
A.~{Sutherland}.
\newblock A database of nonhyperelliptic genus 3 curves over q.
\newblock Available at
  \begin{verbatim}https://math.mit.edu/~drew/gce_genus3_nonhyperelliptic.txt\end{verbatim}.

\bibitem{Weil}
A.~{Weil}.
\newblock {L'arithm\'etique sur les courbes alg\'ebriques.}
\newblock {\em {Acta Math.}}, 52:281--315, 1929.

\end{thebibliography}

\end{document}